\documentclass{amsart}
\usepackage{amssymb}
\usepackage{amscd}

\title{Double shuffle relation for associators}

\author{Hidekazu Furusho}

\address{Graduate School of Mathematics, Nagoya University, 
Furo-cho, Chikusa-ku, Nagoya, 464-8602, Japan}
\email{furusho@math.nagoya-u.ac.jp}

\newcommand{\Exp}{\text{Exp}}
\newtheorem{thm}{Theorem}[section]
\newtheorem{lem}[thm]{Lemma}
\newtheorem{cor}[thm]{Corollary}
\newtheorem{prop}[thm]{Proposition}  

\theoremstyle{remark}
\newtheorem{ack}{Acknowledgments}        

\theoremstyle{definition}

\newtheorem{rem}[thm]{Remark}

\newtheorem{eg}[thm]{Examples}       
\newtheorem{pf}{Proof}

\numberwithin{equation}{section}
\begin{document}
\bibliographystyle{amsalpha+}
\maketitle
\begin{abstract}
It is proved that Drinfel'd's
pentagon equation implies the generalized double shuffle relation.
As a corollary, an embedding
from the Grothendieck-Teichm\"{u}ller group $GRT_1$
into Racinet's double shuffle group $DMR_0$ is obtained,
which settles the project of Deligne-Terasoma.
It is also proved that the gamma factorization formula follows from
the generalized double shuffle relation.
\end{abstract}

\tableofcontents
\setcounter{section}{-1}
\section{Introduction}
This paper shows that Drinfel'd's pentagon equation \cite{Dr} 
implies the generalized double shuffle relation.
As a corollary, we obtain an embedding from 
the Grothendieck-Teichm\"{u}ller (pro-unipotent) group $GRT_1$ (loc.cit)
to Racinet's double shuffle (pro-unipotent) group $DMR_0$ (\cite{R}).
This realizes the project of Deligne-Terasoma \cite{DT} where a different approach is indicated.
Their arguments concern multiplicative convolutions
whereas our methods are based on a bar construction calculus.
We also prove that the gamma factorization formula follows from
the generalized double shuffle relation.
It extends the result in \cite{DT, I} where they show that 
the GT-relations imply the gamma factorization.

Multiple zeta values $\zeta(k_1,\cdots,k_m)$
are the real numbers defined by the following series
$$
\zeta(k_1,\cdots,k_m):=\sum_{0<n_1<\cdots<n_m}\frac{1}{n_1^{k_1}\cdots n_m^{k_m}}
$$
for $m$, $k_1$,\dots, $k_m\in {\bf N} (=\bf{Z}_{>0})$.
This converges if and only if $k_m>1$.
They were studied (allegedly) firstly by Euler \cite{E} for $m=1,2$.
Several types of relations among multiple zeta values have been discussed.
In this paper we focus on two types of relations,
GT-relations and generalized double shuffle relations.
Both of them are described in terms of 
the Drinfel'd associator \cite{Dr}
$$
\varPhi_{KZ}(X_0,X_1)=1+
\sum
(-1)^m\zeta(k_1,\cdots,k_m)X_0^{k_m-1}X_1\cdots X_0^{k_1-1}X_1+
\text{(regularized terms)}
$$
which is a non-commutative formal power series in two variables $X_0$ and $X_1$.
Its coefficients including regularized terms
are explicitly calculated to be linear combinations of 
multiple zeta values in \cite{F03} proposition 3.2.3
by Le-Murakami's method \cite{LM}.
The Drinfel'd associator is introduced as the connection matrix of the Knizhnik-Zamolodchikov equation in \cite{Dr}.

The {\it GT-relations} are a kind of geometric relation.
They consist of 
one pentagon equation \eqref{pentagon} 
and two hexagon equations \eqref{hexagon},\eqref{hexagon-b} (see below)
for group-like (cf. \S \ref{relations}) series.
An associator means a group-like series $\varphi$
satisfying \eqref{pentagon} and for which there exists
$\mu$ such that the pair 
$(\mu, \varphi)$ satisfies
\eqref{hexagon} and \eqref{hexagon-b}.
It is shown in \cite{Dr} that $\varPhi_{KZ}$
satisfies GT-relations (with $\mu=2\pi i$)
by using symmetry of the KZ-system
on configuration spaces.
The following is our previous theorem in \cite{F08}.

\begin{thm}[\cite{F08}theorem 1]
For any group-like series $\varphi$
satisfying \eqref{pentagon} there always exists (unique up to signature)
$\mu\in\bar k$ 
such that \eqref{hexagon} and \eqref{hexagon-b} (see below)
hold for $(\mu,\varphi)$.
\end{thm}

In contrast, the {\it generalized double shuffle relation}
is a kind of combinatorial relation.
It arises from two ways of expressing multiple zeta values
as iterated integrals and as power series.
There are several formulations of the relations (see \cite{IKZ, R}).
In particular, they are formulated as \eqref{double shuffle} (see below)
for $\varphi=\varPhi_{KZ}$ in \cite{R}.

\begin{thm}\label{main}
Let $\varphi$ be a non-commutative formal power series in two variables
which is group-like.
Suppose that $\varphi$ satisfies Drinfel'd's pentagon equation \eqref{pentagon}.
Then it also satisfies the generalized double shuffle relation \eqref{double shuffle} (see below).
\end{thm}

We note that a similar result is announced by Terasoma in \cite{T}.
The essential part of our proof is to use the series shuffle formula 
\eqref{series shuffle for Li}(see below),
a functional relation among complex multiple polylogarithms 
\eqref{one variable MPL} and \eqref{two variables MPL}.
This induces the series shuffle formula 
\eqref{series shuffle for bar}
for the corresponding elements 
in the bar construction of the moduli space ${\mathcal M}_{0,5}$.
We evaluate each term of \eqref{series shuffle for bar}
at the product of the last two terms of \eqref{5-cycle}
in \S\ref{auxiliary}
and conclude the series shuffle formula \eqref{series shuffle for l}
for above $\varphi$.

The {\it Grothendieck-Teichm\"{u}ller group} $GRT_1$
is a pro-unipotent group introduced by Drinfel'd \cite{Dr}
which is closely related to Grothendieck's philosophy of
Teichm\"{u}ller-Lego in \cite{Gr}.
Its set of $k$-valued ($k$: a field with characteristic $0$) points is
defined to be the set of associators $\varphi$ with $\mu=0$.
Its multiplication is given by
\begin{equation}\label{multiplication}
\varphi_1\circ\varphi_2:=\varphi_1(\varphi_2 X_0\varphi_2^{-1},X_1)\cdot\varphi_2
=\varphi_2\cdot\varphi_1(X_0,\varphi^{-1}_2X_1\varphi_2).
\end{equation}
In contrast, the {\it double shuffle group} $DMR_0$
is a pro-unipotent group introduced by Racinet \cite{R}.
Its set of $k$-valued points
consists of group-like series $\varphi$ which satisfy \eqref{double shuffle}
\footnote{
For our convenience,
we change some signatures in the original definition (\cite{R} definition 3.2.1.))
}
and $c_{X_0}(\varphi)=c_{X_1}(\varphi)=c_{X_0 X_1}(\varphi)=0$.
(For a monic monomial $W$, 
$c_W(\varphi)$ is the coefficient of $W$ in $\varphi$.)
Its multiplication
\footnote{
Again for our convenience,
we change the order of multiplication in \cite{R} (3.1.2.1.)
}
is given by the equation \eqref{multiplication}.
It is conjectured that both groups are isomorphic
to the unipotent part of the motivic Galois group of $\bf Z$;
the Galois group of unramified mixed Tate motives
(as explained in \cite{A}).
The following are direct corollaries of our theorem \ref{main}
since the equations \eqref{hexagon} and \eqref{hexagon-b} for $(\mu,\varphi)$
imply $c_{X_0X_1}(\varphi)=\frac{\mu^2}{24}$.

\begin{cor}
The Grothendieck-Teichm\"{u}ller group $GRT_1$
is embedded in the double shuffle group $DMR_0$
as pro-algebraic groups.
\end{cor}

Considering their associated Lie algebras,
we get an embedding from the Grothendieck-Teichm\"{u}ller Lie algebra
$\frak{grt}_1$ \cite{Dr} 
into the double shuffle Lie algebra
$\frak{dmr}_0$ \cite{R}.

\begin{cor}
For $\mu\in k^\times$,
the Grothendieck-Teichm\"{u}ller torsor $M_\mu$
is embedded in the double shuffle torsor $DMR_\mu$ as pro-torsors.
\end{cor}

Here $M_\mu$ is the right $GRT_1$-torsor in \cite{Dr} 
whose action is given by the equation \eqref{multiplication} and
whose set of $k$-valued points is defined to be
the collection of associators $\varphi$ with $\mu$ and 
$DMR_\mu$ is the right $DMR_0$-torsor in \cite{R}
whose action is given by the equation \eqref{multiplication} and
whose set of $k$-valued points is defined to
the collection of group-like series $\varphi$ which satisfy \eqref{double shuffle},
$c_{X_0}(\varphi)=c_{X_1}(\varphi)=0$ and
$c_{X_0 X_1}(\varphi)=\frac{\mu^2}{24}$ respectively.
We note that $\varPhi_{KZ}$ gives an element of $M_\mu(\bf C)$ and $DMR_\mu(\bf C)$
with $\mu=2\pi i$.

Let $\varphi\in k\langle\langle X_0,X_1\rangle\rangle$
be a non-commutative formal power series in two variables
which is group-like with $c_{X_0}(\varphi)=c_{X_1}(\varphi)=0$.
It is uniquely expressed as
$\varphi=1+\varphi_{X_0} X_0+\varphi_{X_1} X_1$
with $\varphi_{X_0}$, $\varphi_{X_1}\in k\langle\langle X_0,X_1\rangle\rangle$.
The meta-abelian quotient $B_\varphi(x_0,x_1)$ of $\varphi$ is defined to be
$(1+\varphi_{X_1} X_1)^{\text{ab}}$
where $h\mapsto h^{\text{ab}}$ is the abelianization map
$k\langle\langle X_0,X_1\rangle\rangle\to k[[x_0,x_1]]$.
The following is our second main theorem.

\begin{thm}\label{gamma factorization}
Let $\varphi$ 
be a non-commutative formal power series in two variables
which is group-like with $c_{X_0}(\varphi)=c_{X_1}(\varphi)=0$.
Suppose that it satisfies the generalized double shuffle relation \eqref{double shuffle}.
Then its meta-abelian quotient is gamma-factorisable, i.e.
there exists a unique series $\Gamma_\varphi(s)$ in
$1+s^2 k[[s]]$ such that
\begin{equation}\label{gamma factorization formula}
B_\varphi(x_0,x_1)=\frac{\Gamma_\varphi(x_0)\Gamma_\varphi(x_1)}
{\Gamma_\varphi(x_0+x_1)}.
\end{equation}
The gamma element $\Gamma_\varphi$ gives the correction term
$\varphi_{\text{corr}}$ of the series shuffle regularization \eqref{correction}
by $\varphi_{\text{corr}}=\Gamma_\varphi(-Y_1)^{-1}$.
\end{thm}

This theorem extends the results in \cite{DT, I} which show that 
for any group-like series satisfying \eqref{pentagon}, \eqref{hexagon} and \eqref{hexagon-b}
its meta-abelian quotient is gamma factorisable.
This result might be a step to relate  $DMR$ 
with the set $\text{Sol}KV$
of solutions of the Kashiwara-Vergne equations
which is defined by a certain tangential automorphism condition and
a coboundary Jacobian condition \cite{AT, AET}
since the equation \eqref{gamma factorization formula} is regarded as a consequence of
the latter condition (cf. \cite{AET}\S 2.1.)
It is calculated in \cite{Dr} that
$\Gamma_\varphi(s)=\exp\{{\sum_{n=2}^\infty}\frac{\zeta(n)}{n}s^n\}
=e^{-\gamma s}\Gamma(1-s)$
($\gamma$: Euler constant, $\Gamma(s)$: the classical gamma function)
in particular case of $\varphi=\varPhi_{KZ}$.

\S\ref{relations} is a review of the GT-relations and
the generalized double shuffle relation.
\S\ref{bar} is to recall bar constructions 
which is a main tool of the proof of theorem \ref{main} in 
\S\ref{proof of main theorem}.
Auxiliary lemmas which are necessary to the proof
are shown in \S\ref{auxiliary}.
Theorem \ref{gamma factorization} is proved 
in \S\ref{proof of gamma factorization}.
Appendix \ref{appendix} is  a brief review of the essential part of
the proof of Racinet's
theorem that $DMR_0$ forms a group. 

\section{GT-relations and generalized double shuffle relation}
\label{relations}
Let $k$ be a field with characteristic $0$.
Let $U\frak F_2=k\langle\langle X_0,X_1\rangle\rangle$ be a non-commutative formal power series ring
in two variables $X_0$ and $X_1$.
An element $\varphi=\varphi(X_0,X_1)$ is called {\it group-like} if it satisfies
$
\Delta (\varphi)=\varphi\widehat\otimes \varphi
$
with 
$\Delta(X_0)=X_0\otimes 1+1\otimes X_0$ and
$\Delta(X_1)=X_1\otimes 1+1\otimes X_1$
and its constant term is $1$.
Here, $\widehat\otimes$ means the completed tensor product.
For any $k$-algebra homomorphism $\iota:U\frak F_2\to S$
the image $\iota(\varphi)\in S$ is denoted 
by $\varphi(\iota(X_0),\iota(X_1))$.
Let $\frak a_4$ be the completion (with respect to the natural grading)
of the pure braid Lie algebra with 4-strings;
the Lie algebra over $k$ with generators
$t_{ij}$ ($1\leqslant i,j \leqslant 4$) 
and defining relations $t_{ii}=0$,
$t_{ij}=t_{ji}$, $[t_{ij},t_{ik}+t_{jk}]=0$ ($i$,$j$,$k$: all distinct)
and $[t_{ij},t_{kl}]=0$ ($i$,$j$,$k$,$l$: all distinct).
Let $\varphi=\varphi(X_0,X_1)$ be a 
group-like element of $U\frak F_2$ and $\mu\in k$.
The {\it GT-relations} for $(\mu,\varphi)$ consist of one pentagon equation
\begin{equation}\label{pentagon}
\varphi(t_{12},t_{23}+t_{24})
\varphi(t_{13}+t_{23},t_{34})=
\varphi(t_{23},t_{34})
\varphi(t_{12}+t_{13},t_{24}+t_{34})
\varphi(t_{12},t_{23})
\end{equation}
and two hexagon equations
\begin{equation}\label{hexagon}
\exp\{\frac{\mu (t_{13}+t_{23})}{2}\}=
\varphi(t_{13},t_{12})\exp\{\frac{\mu t_{13}}{2}\}
\varphi(t_{13},t_{23})^{-1}
\exp\{\frac{\mu t_{23}}{2}\} \varphi(t_{12},t_{23}),
\end{equation}
\begin{equation}\label{hexagon-b}
\exp\{\frac{\mu (t_{12}+t_{13})}{2}\}=
\varphi(t_{23},t_{13})^{-1}\exp\{\frac{\mu t_{13}}{2}\}
\varphi(t_{12},t_{13})
\exp\{\frac{\mu t_{12}}{2}\} \varphi(t_{12},t_{23})^{-1}.
\end{equation}

Let $\pi_Y:k\langle\langle X_0, X_1\rangle\rangle\to
k\langle\langle Y_1,Y_2,\dots\rangle\rangle$ be 
the $k$-linear map between non-commutative formal power series rings
that sends all the words ending in $X_0$ to zero and the
word $X_0^{n_m-1}X_1\cdots X_0^{n_1-1}X_1$ ($n_1,\dots,n_m\in\bold N$)
to $(-1)^mY_{n_m}\cdots Y_{n_1}$.
Define the coproduct $\Delta_*$ on $k\langle\langle
Y_1,Y_2,\dots\rangle\rangle$ by
$\Delta_* Y_n=\sum_{i=0}^n Y_i\otimes Y_{n-i}$ with $Y_0:=1$.
For $\varphi=\sum_{W:\text{word}} c_W(\varphi) W\in
k\langle\langle X_0,X_1\rangle\rangle$, 
define the series shuffle regularization
$\varphi_*=\varphi_{\text{corr}}\cdot\pi_Y(\varphi)$
with the correction term
\begin{equation}\label{correction}
\varphi_{\text{corr}}=\exp\left(\sum_{n=1}^{\infty}
\frac{(-1)^n}{n}c_{X_0^{n-1}X_1}(\varphi)Y_1^n\right).
\end{equation}
For a group-like series $\varphi\in U\frak F_2$
the {\it generalised double shuffle relation} means the equality
\begin{equation}\label{double shuffle}
\Delta_*(\varphi_*)=\varphi_*\widehat\otimes \varphi_*.
\end{equation}

Let ${\frak P_5}$ stand for the completion 
(with respect to the natural grading) of 
the pure sphere braid Lie algebra
with 5 strings;
the Lie algebra over $k$ 
generated by $X_{ij}$ ($1\leqslant i,j\leqslant 5$) with relations
$X_{ii}=0$, $X_{ij}=X_{ji}$, $\sum_{j=1}^5 X_{ij}=0$
($1\leqslant i,j\leqslant 5$) and
$[X_{ij},X_{kl}]=0$ if $\{i,j\}\cap\{k,l\}=\emptyset$.
Denote its universal enveloping algebra by $U{\frak P_5}$.
(Note: $X_{45}=X_{12}+X_{13}+X_{23}$, $X_{51}=X_{23}+X_{24}+X_{34}$.)
There is a surjection $\tau:\frak a_4\to \frak P_5$
sending $t_{ij}$ to $X_{ij}$ ($1\leqslant i,j\leqslant 4$).
Its kernel is the center of $\frak a_4$ generated by 
$\sum_{1\leqslant i,j\leqslant 4}t_{ij}$.
By \cite{F08} lemma 5, theorem \ref{main} is reduced to the following.

\begin{thm}\label{almost main}
Let $\varphi$ be a group-like element of $U\frak F_2$
with $c_{X_0}(\varphi)=c_{X_1}(\varphi)=0$.
Suppose that $\varphi$ satisfies the 5-cycle relation in $U{\frak P_5}$
\begin{equation}\label{5-cycle}
\varphi(X_{34},X_{45})
\varphi(X_{51},X_{12})
\varphi(X_{23},X_{34})
\varphi(X_{45},X_{51})
\varphi(X_{12},X_{23})
=1.
\end{equation}
Then it also satisfies the generalized double shuffle relation,
i.e.
$\Delta_*(\varphi_*)=\varphi_*\widehat\otimes \varphi_*$.
\end{thm}

\section{Bar constructions}\label{bar}
In this section we review the notion of bar construction and
multiple polylogarithm functions
which are essential to prove our main theorem.

Let $\mathcal M_{0,4}$ be the moduli space
$\{(x_1,\cdots,x_4)\in (\bold P_k^1)^4|x_i\neq x_j (i\neq j)\}/ PGL_2(k)$ 
of 4 different points in $\bold P^1$.
It is identified with
$\{z\in\bold P^1 | z\neq 0,1,\infty\}$
by sending $[(0,z,1,\infty)]$ to $z$.
Let $\mathcal M_{0,5}$ be the moduli space
$\{(x_1,\cdots,x_5)\in (\bold P_k^1)^5|x_i\neq x_j (i\neq j)\}/ PGL_2(k)$
of 5 different points in $\bold P^1$.
It is identified with
$\{(x,y)\in\bold G_m^2 | x\neq 1,y\neq 1,xy\neq 1\}$
by sending $[(0,xy,y,1,\infty)]$ to $(x,y)$.

For $\mathcal M=\mathcal M_{0,4}/k$ or $\mathcal M_{0,5}/k$,
we consider Brown's variant $V(\mathcal M)$ \cite{B} of
Chen's reduced bar construction \cite{C}.
This is a graded Hopf algebra
$V(\mathcal M)=\oplus_{m=0}^\infty V_m$
($\subset TV_1=\oplus_{m=0}^\infty V_1^{\otimes m}$)
over $k$.
Here  $V_0=k$, $V_1=H^1_{DR}(\mathcal M)$ and
$V_m$ is the totality of linear combinations (finite sums)
$\sum_{I=(i_m,\cdots,i_1)}c_I[\omega_{i_m}|\cdots|\omega_{i_1}]
\in V_1^{\otimes m}$
($c_I\in k$, $\omega_{i_j}\in V_1$,
$[\omega_{i_m}|\cdots|\omega_{i_1}]:=\omega_{i_m}\otimes\cdots\otimes\omega_{i_1}$)
satisfying the integrability condition
\begin{equation}\label{integrability}
\sum_{I=(i_m,\cdots,i_1)}c_I[\omega_{i_m}|\omega_{i_{m-1}}|\cdots
|\omega_{i_{j+1}}\wedge\omega_{i_{j}}|\cdots|\omega_{i_1}]=0
\end{equation}
in $V_1^{\otimes m-j-1}\otimes H^2_{DR}(\mathcal M)\otimes V_1^{\otimes j-1}$
for all $j$ ($1\leqslant j<m$). 

For $\mathcal M=\mathcal M_{0,4}$,
$V(\mathcal M_{0,4})$ is generated by $\omega_0=d\log(z)$ and $\omega_1=d\log(z-1)$
and $\omega_0\wedge\omega_1=0$.
We identify $V({\mathcal M}_{0,4})$ with the graded $k$-linear dual of
$U\frak F_2$  by
$\Exp\Omega_4=\Exp\Omega_4(x;X_0,X_1):=\sum X_{i_m}\cdots X_{i_1}\otimes
[\omega_{i_m}|\cdots |\omega_{i_1}]\in
U\frak F_2\widehat\otimes V({\mathcal M}_{0,4})$.
Here the sum is taken over $m\geqslant 0$ and $i_1,\cdots,i_m\in\{0,1\}$.
It is easy to see that the identification is compatible
with Hopf algebra structures.
We note that the product $l_1\cdot l_2\in V(\mathcal M_{0,4})$ 
for $l_1$, $l_2\in V(\mathcal M_{0,4})$ is given by
$l_1\cdot l_2(f):=\sum_i l_1(f^{(i)}_1) l_2(f^{(i)}_2)$ for
$f\in U\frak F_2$ with $\Delta(f)=\sum_i f_1^{(i)}\otimes f_2^{(i)}$.
Occasionally we regard $V({\mathcal M}_{0,4})$ as the regular function ring
of $F_2(k)=\{g\in U\frak F_2| g:\text{group-like}\}
=\{g\in U\frak F_2| g(0)=1, \Delta(g)=g\otimes g\}$.

For $\mathcal M=\mathcal M_{0,5}$,
$V(\mathcal M_{0,5})$ is generated by $\omega_{k,l}=d\log (x_k-x_l)$ ($1\leqslant k,l\leqslant 5$)
with relations $\omega_{ii}=0$, $\omega_{ij}=\omega_{ji}$, $\sum_{j=1}^5 \omega_{ij}=0$
($1\leqslant i,j\leqslant 5$) and
$\omega_{ij}\wedge\omega_{kl}=0$ if $\{i,j\}\cap\{k,l\}=\emptyset$.
By definition $V({\mathcal M}_{0,5})$ is naturally identified 
with the graded $k$-linear dual of $U\frak P_5$.
The identification is induced from 
$$
\Exp\Omega_5:=\sum X_{J_m}\cdots X_{J_1}\otimes
[\omega_{J_m}|\cdots |\omega_{J_1}]\in
U\frak P_5\widehat\otimes TV_1
$$
where the sum is taken over $m\geqslant 0$ and $J_1,\cdots,J_m\in
\{(k,l)|1\leqslant k<l\leqslant 5\}$.

\begin{rem}
Here every monomial in $U\frak P_5$ appears in the left-hand tensor factor of
$\Exp\Omega_5$, but when these monomials are gathered in terms of a linear basis of
$U\frak P_5$, the right-hand tensor factors automatically gather into
linear combinations which forms a basis of $V({\mathcal M}_{0,5})$,
and all of which satisfies \eqref{integrability}.
Hence $\Exp\Omega_5$ lies on $U\frak P_5\widehat\otimes V({\mathcal M}_{0,5})$.
\end{rem}

Especially the identification between degree 1 terms is given by
$$
\Omega_5=\sum_{1\leqslant k<l\leqslant 5} X_{kl}d\log(x_k-x_l)
\in \frak P_5\otimes H^1_{DR}({\mathcal M}_{0,5}).
$$
In terms of the coordinate $(x,y)$,
\begin{align*}
\Omega_5
&=X_{12}d\log (xy)+X_{13}d\log y+X_{23}d\log y(1-x)+X_{24}d\log (1-xy)+X_{34}d\log (1-y)\\
&=X_{12}d\log x+X_{23}d\log (1-x)+(X_{12}+X_{13}+X_{23})d\log y+X_{34}d\log (1-y)\\
&\qquad\qquad +X_{24}d\log (1-xy)\\
&=X_{12}\frac{dx}{x}+X_{23}\frac{dx}{x-1}+X_{45}\frac{dy}{y}+X_{34}\frac{dy}{y-1}+X_{24}\frac{ydx+xdy}{xy-1}.
\end{align*}
It is easy to see that the identification is compatible
with Hopf algebra structures.
We note again that
the product $l_1\cdot l_2\in V(\mathcal M_{0,5})$
for $l_1$, $l_2\in V(\mathcal M_{0,5})$ is given by
$l_1\cdot l_2(f):=\sum_i l_1(f^{(i)}_1) l_2(f^{(i)}_2)$ for
$f\in U\frak P_5$ with $\Delta(f)=\sum_i f_1^{(i)}\otimes f_2^{(i)}$
($\Delta$: the coproduct of $U\frak P_5$).
Occasionally we also regard $V({\mathcal M}_{0,5})$ as the regular function ring
of $P_5(k)=\{g\in U\frak P_5| g:\text{group-like}\}$.

For the moment assume that $k$ is a subfield of $\bold C$.
We have an embedding (called a realisation in  \cite{B}\S 1.2, \S 3.6)
$\rho:V(\mathcal M)\hookrightarrow I_o(\mathcal M)$
as algebra over $k$ which sends 
$\sum_{I=(i_m,\cdots,i_1)}c_I[\omega_{i_m}|\cdots|\omega_{i_1}]$ 
($c_I\in k$) to
$\sum_{I}c_I\text{It}\int_o\omega_{i_m}\circ\cdots\circ\omega_{i_1}$.
Here 
$\sum_{I}c_I\text{It}\int_o\omega_{i_m}\circ\cdots\circ\omega_{i_1}$
means the iterated integral defined by
\begin{equation}\label{iterated integral map}
\sum_{I}c_I\int_{0<t_1< \cdots <t_{m-1}<t_m<1}
\omega_{i_m}({\gamma(t_m)})\cdot
\omega_{i_{m-1}}({\gamma(t_{m-1})})\cdot
\cdots\omega_{i_1}({\gamma(t_1)})
\end{equation}
for all analytic paths $\gamma: (0,1)\to \mathcal M(\bold C)$
starting from the tangential basepoint $o$
(defined by $\frac{d}{dz}$ for $\mathcal M=\mathcal M_{0,4}$ and
defined by $\frac{d}{dx}$ and $\frac{d}{dy}$ for $\mathcal M=\mathcal M_{0,5}$)
at the origin in $\mathcal M$ (for its treatment see also \cite{De}\S 15)
and $I_o(\mathcal M)$
\footnote{
In \cite{B} it is denoted by $L_o({\mathcal M})$.
}
denotes 
the ${\mathcal O}^\text{an}_{\mathcal M}$-module generated by all such 
homotopy invariant iterated integrals with $m\geqslant 1$ and holomorphic 1-forms
$\omega_{i_1},\dots,\omega_{i_m}\in \Omega^1(\mathcal M)$.

For ${\bf a}=(a_1,\cdots,a_k)\in\bold Z^k_{>0}$,
its weight and its depth are defined to be $wt({\bf a})=a_1+\cdots+a_k$ and
$dp({\bf a})=k$ respectively. 
Put $z\in\bf C$ with $|z|<1$.
Consider the following complex function which is 
called the {\it one variable multiple polylogarithm}
\begin{equation}\label{one variable MPL}
Li_{\bold a}(z):=\underset{0<m_1<\cdots<m_k}{\sum}
\frac{z^{m_k}}{m_1^{a_1}\cdots m_k^{a_k}}.
\end{equation}
It satisfies the following differential equation
$$
\frac{d}{dz}Li_{\bold a}(z)=
\begin{cases}
\frac{1}{z}Li_{(a_1,\cdots,a_{k-1},a_k-1)}(z)
&\text{if } a_k\neq 1,\\
\frac{1}{1-z}Li_{(a_1,\cdots,a_{k-1})}(z)
&\text{if } a_k=1,k\neq 1, \\
\frac{1}{1-z}
&\text{if } a_k=1,k=1. \\
\end{cases}\notag
$$
It gives an iterated integral starting from $o$,
which lies on $I_o(\mathcal M_{0,4})$.
Actually it corresponds to an element of
$V(\mathcal M_{0,4})$ denoted by $l_{\bf a}$.
It is expressed as
\begin{equation}\label{Li-expression}
l_{\bf a}=(-1)^k[\underbrace{\omega_0|\cdots|\omega_0}_{a_k-1}|\omega_1|
\underbrace{\omega_0|\cdots|\omega_0}_{a_{k-1}-1}|\omega_1|
\omega_0|\cdots\cdots|\omega_1|
\underbrace{\omega_0|\cdots|\omega_0}_{a_1-1}|\omega_1]
\end{equation}
and is calculated by
$l_{\bf a}(\varphi)=(-1)^k c_{X_0^{a_k-1}X_1X_0^{a_{k-1}-1}X_1\cdots X_0^{a_1-1}X_1}(\varphi)$
for a series $\varphi=\sum_{W:\text{word}} c_W(\varphi) W$.

For ${\bf a}=(a_1,\cdots,a_k)\in\bold Z^k_{>0}$,
${\bf b}=(b_1,\cdots,b_l)\in\bold Z^l_{>0}$
and $x,y\in\bf C$ with $|x|<1$ and $|y|<1$,
consider the following complex function which is called
the {\it two variables multiple polylogarithm}
\begin{equation}\label{two variables MPL}
Li_{\bold a,\bold b}(x,y):=
\underset{<n_1<\cdots<n_l}{\underset{0<m_1<\cdots<m_k}{\sum}}
\frac{\qquad x^{m_k} \qquad y^{n_l}}
{m_1^{a_1}\cdots m_k^{a_k}n_1^{b_1}\cdots n_l^{b_l}}.
\end{equation}
It satisfies the following differential equations in \cite{BF}\S 5

\begin{align}\label{differential equation}
&\frac{d}{dx}Li_{\bold a,\bold b}(x,y)=
\begin{cases}
\frac{1}{x}Li_{(a_1,\cdots,a_{k-1},a_k-1),\bold b}(x,y)
\qquad\qquad\qquad\qquad \ \ 
\text{if } a_k\neq 1,\\
\frac{1}{1-x}Li_{(a_1,\cdots,a_{k-1}),\bold b}(x,y)
-\left(\frac{1}{x}+\frac{1}{1-x}\right)
Li_{(a_1,\cdots,a_{k-1},b_1),(b_2,\cdots,b_l)}(x,y)\\
\qquad\qquad\qquad\qquad\qquad\qquad\qquad\qquad\qquad\qquad
\text{if } a_k=1,k\neq 1, l\neq 1,\\
\frac{1}{1-x}Li_{\bold b}(y)
-\left(\frac{1}{x}+\frac{1}{1-x}\right)
Li_{(b_1),(b_2,\cdots,b_l)}(x,y)
\text{   if } a_k=1,k= 1, l\neq 1,\\
\frac{1}{1-x}Li_{(a_1,\cdots,a_{k-1}),\bold b}(x,y)
-\left(\frac{1}{x}+\frac{1}{1-x}\right)
Li_{(a_1,\cdots,a_{k-1},b_1)}(xy)\\
\qquad\qquad\qquad\qquad\qquad\qquad\qquad\qquad\qquad\qquad
\text{if } a_k=1,k\neq 1, l=1,\\
\frac{1}{1-x}Li_{\bold b}(y)
-\left(\frac{1}{x}+\frac{1}{1-x}\right)
Li_{\bold b}(xy)
\qquad\qquad\quad \ \ 
\text{if } a_k=1,k=1, l=1,\\
\end{cases}
\\
&\frac{d}{dy}Li_{\bold a,\bold b}(x,y)=
\begin{cases}
\frac{1}{y}Li_{\bold a,(b_1,\cdots,b_{l-1},b_l-1)}(x,y)
&\text{if } b_l\neq 1,\\
\frac{1}{1-y}Li_{\bold a,(b_1,\cdots,b_{l-1})}(x,y)
&\text{if } b_l=1, l\neq 1,\\
\frac{1}{1-y}Li_{\bold a}(xy)
&\text{if } b_l=1,l=1.\\
\end{cases}\notag 
\end{align}
By analytic continuation, the functions 
$Li_{\bf a,\bf b}(x,y)$, $Li_{\bf a,\bf b}(y,x)$,
$Li_{\bf a}(x)$, $Li_{\bf a}(y)$ and $Li_{\bf a}(xy)$
give iterated integrals starting from $o$,
which lie on $I_o(\mathcal M_{0,5})$.
They correspond to elements of $V(\mathcal M_{0,5})$ by the map $\rho$
denoted by
$l_{\bf a,\bf b}^{x,y}$, $l_{\bf a,\bf b}^{y,x}$,
$l_{\bf a}^x$, $l_{\bf a}^y$  and $l_{\bf a}^{xy}$ respectively.
Note that they are expressed as
\begin{equation}\label{expression}
\sum_{I=(i_m,\cdots,i_1)}c_I[\omega_{i_m}|\cdots|\omega_{i_1}]
\end{equation}
for some $m\in\bf N$
with $c_I\in \bold Q$ and $\omega_{i_j}\in \{\frac{dx}{x},\frac{dx}{1-x},\frac{dy}{y},\frac{dy}{1-y},\frac{xdy+ydx}{1-xy}\}$.
It is easy to see that they  indeed lie on $V(\mathcal M_{0,5})$,
i.e. they satisfy the integrability condition \eqref{integrability},
by induction on weights:
Suppose that $l_{\bf a,\bf b}^{x,y}$ is expressed as above.
Then by our induction assumption, \eqref{integrability} holds for
$1\leqslant j\leqslant m-2$.
By the integrability of the complex analytic function $Li_{\bf a,\bf b}(x,y)$,
we have $d \circ d Li_{\bf a,\bf b}(x,y)=0$.
This implies  \eqref{integrability} for $j=m-1$.
The same arguments also work for $l_{\bf a,\bf b}^{y,x}$,
$l_{\bf a}^x$, $l_{\bf a}^y$  and $l_{\bf a}^{xy}$. 

\begin{eg}\label{first example}
Put 
$\alpha_0=\frac{dx}{x}$, $\alpha_1=\frac{dx}{1-x}$,
$\beta_0=\frac{dy}{y}$, $\beta_1=\frac{dy}{1-y}$ and
$\gamma=\frac{xdy+ydx}{1-xy}$.
The function $Li_{2,1}(x,y)$ corresponds to
$l^{x,y}_{2,1}=
[\beta_1|\alpha_0|\gamma]
+[\beta_1|\beta_0|\gamma]
+[\alpha_0|\beta_1|\gamma]
+[\alpha_0|\alpha_1|\beta_1]
-[\alpha_0|\alpha_0|\gamma]
-[\alpha_0|\alpha_1|\gamma]
$.
By a direct computation it can be checked that $l^{x,y}_{2,1}$
lies on $V(\mathcal M_{0,5})$.
The analytic continuation $\rho(l^{x,y}_{2,1})$ 
of $Li_{2,1}(x,y)$ is calculated by \eqref{iterated integral map}.
In particular when we take any path $\gamma(t)$ starting from $o$ to $(x,y)$ in the open unit disk
of $\mathcal M_{0,5}(\bf C)$ we get the expression \eqref{two variables MPL}
of $Li_{2,1}(x,y)$ for $|x|,|y|<1$.
\end{eg}

\section{Proof of theorem \ref{main}}\label{proof of main theorem}
This section gives a proof of theorem \ref{main}.

Suppose that $\varphi$ is an element as in theorem \ref{almost main}.
Recall that multiple polylogarithms satisfy the analytic identity,
the series shuffle formula in
$I_o(\mathcal M_{0,5})$
\begin{equation}\label{series shuffle for Li}
Li_{\bold a}(x)\cdot Li_{\bold b}(y)
={\underset{\sigma\in Sh^{\leqslant}(k,l)}{\sum}}
Li_{\sigma(\bold a,\bold b)}(\sigma(x,y)).
\end{equation}
Here
$
Sh^{\leqslant}(k,l):={\cup}^\infty_{N=1}\{
\sigma:\{1,\cdots,k+l\}\to\{1,\cdots,N \}| \sigma{\text{ is onto}},
\sigma(1)<\cdots<\sigma(k), \sigma(k+1)<\cdots<\sigma(k+l)
\},
$
$\sigma(\bold a,\bold b):=((c_1,\cdots, c_j),(c_{j+1},\cdots,c_{N}))$ 
with $\{j,N\}=\{\sigma(k),\sigma(k+l)\}$,
$$
c_i=
\begin{cases}
a_s+b_{t-k}  &\text{if } \sigma^{-1}(i)=\{s,t\} \text{ with } s<t , \\
a_s     &\text{if } \sigma^{-1}(i)=\{s\}  \quad \text{with } s\leqslant k,\\
b_{s-k} &\text{if } \sigma^{-1}(i)=\{s\}  \quad \text{with } s> k,\\
\end{cases}
\text{ and }
\sigma(x,y)=
\begin{cases}
xy    &\text{if } \sigma^{-1}(N)=k,k+l,\\
(x,y) &\text{if } \sigma^{-1}(N)=k+l,\\
(y,x) &\text{if } \sigma^{-1}(N)=k.
\end{cases}
$$
Since $\rho$ is an embedding of algebras,
the above analytic identity immediately implies the algebraic identity,
the series shuffle formula in $V(\mathcal M_{0,5})$
\begin{equation}\label{series shuffle for bar}
l_{\bold a}^x\cdot l_{\bold b}^y
={\underset{\sigma\in Sh^{\leqslant}(k,l)}{\sum}}
l_{\sigma(\bold a,\bold b)}^{\sigma(x,y)}.
\end{equation}
\begin{eg}
Under the notation in example \ref{first example},
$l_2^x=[\alpha_0|\alpha_1]$, $l_1^y=[\beta_1]$,
$l_3^{xy}=[\alpha_0+\beta_0|\alpha_0+\beta_0|\gamma]$ and
$l_{1,2}^{y,x}=[\beta_1|\alpha_0|\alpha_1]-[\beta_0|\alpha_0|\gamma]
-[\beta_0|\beta_0|\gamma]-[\beta_1|\alpha_0|\gamma]
-[\beta_1|\beta_0|\gamma]+[\alpha_0|\beta_1|\alpha_1]
-[\alpha_0|\beta_0|\gamma]-[\alpha_0|\beta_1|\gamma]
+[\alpha_0|\alpha_1|\gamma]$.
It can be checked the equation \eqref{series shuffle for bar}
in the case of $\bold a=(2)$ and $\bold b=(1)$ by a direct calculation
(for $l_{2,1}^{x,y}$ see example \ref{first example}.)
\end{eg}
Evaluation of the equation \eqref{series shuffle for bar}
at the group-like element $\varphi_{451}\varphi_{123}$
\footnote{
For simplicity we mean $\varphi_{ijk}$ for
$\varphi(X_{ij},X_{jk})\in U\frak P_5$.}
gives the series shuffle formula
\begin{equation}\label{series shuffle for l}
l_{\bold a}(\varphi)\cdot l_{\bold b}(\varphi)
={\underset{\sigma\in Sh^{\leqslant}(k,l)}{\sum}}
l_{\sigma(\bold a,\bold b)}(\varphi)
\end{equation}
for admissible
\footnote{
An index ${\bf a}=(a_1,\cdots,a_k)$ is called {\it admissible} if $a_k>1$.
}
indices $\bf a$ and $\bf b$
because of lemma \ref{lemma3} and lemma \ref{lemma4} below.

Define the integral regularized value $l^I_{\bf a}(\varphi)$ 
and series regularized value $l^S_{\bf a}(\varphi)$ 
in $k[T]$ ($T$: a parameter which stands for $\log x$)
for all indices $\bf a$ as follows:
$l^I_{\bf a}(\varphi)$ is defined by 
$l^I_{\bf a}(\varphi)=\l_{\bf a}(e^{TX_1}\varphi)$.
(Equivalently $l^I_{\bf a}(\varphi)$ for any index $\bf a$
is uniquely defined in such a way that
the iterated integral shuffle formulae (cf.\cite{BF}\S 7)
remain valid 
for all indices $\bf a$
with $l^I_1(\varphi):=-T$ and
$l^I_{\bf a}(\varphi):=l_{\bf a}(\varphi)$ for all admissible indices $\bf a$.)
Similarly put $l^S_1(\varphi):=-T$ and
put $l^S_{\bf a}(\varphi):=l_{\bf a}(\varphi)$ for all admissible indices $\bf a$.
Then by the results in \cite{H},
$l^S_{\bf a}(\varphi)$ for a non-admissible index $\bf a$ is also uniquely defined
in such a way that the series shuffle formulae \eqref{series shuffle for l} 
remain valid 
for $l^S_{\bf a}(\varphi)$ with all indices $\bf a$.

Let ${\mathbb L}$ be the $k$-linear map from
$k[T]$ to itself defined via the generating function:
\begin{equation}\label{generating function}
{\mathbb L}(\exp Tu)
=\sum_{n=0}^{\infty}{\mathbb L}(T^n)\frac{u^n}{n!}
=\mbox{exp}\left\{-\sum_{n=1}^{\infty}{l^I_n(\varphi)}\frac{u^n}{n}\right\}
\left(=\mbox{exp}\left\{Tu-\sum_{n=1}^{\infty}{l_n(\varphi)}\frac{u^n}{n}\right\}\right).
\end{equation}

\begin{prop}
Let $\varphi$ be an element as in theorem \ref{almost main}.
Then the regularization relation holds, i.e.
$\l^S_{\bf a}(\varphi)={\mathbb L}\bigl(l^I_{\bf a}(\varphi)\bigr)$
for all indices $\bf a$.
\end{prop}

\begin{pf}
We may assume that $\bf a$ is non-admissible
because the proposition is trivial if $\bf a$ is admissible.
When ${\bf a}$ is of the form $(1,1,\cdots,1)$,
the proof is given by the same argument as in \cite{G} lemma 7.9 as follows:
By the series shuffle formulae,
$$
\sum_{k=0}^m (-1)^kl^S_{k+1}(\varphi)\cdot l^S_{\underbrace{1,1,\cdots,1}_{m-k}}(\varphi)
=(m+1)l^S_{\underbrace{1,1,\cdots,1}_{m+1}}(\varphi)
$$
for $m\geqslant 0$. Here we put $l^S_\emptyset(\varphi)=1$.
This means
$$
\sum_{k,l\geqslant 0} (-1)^kl^S_{k+1}(\varphi)\cdot l^S_{\underbrace{1,1,\cdots,1}_{l}}(\varphi)u^{k+l}
=\sum_{m\geqslant 0}(m+1)l^S_{\underbrace{1,1,\cdots,1}_{m+1}}(\varphi)u^m.
$$
Put
$f(u)=\sum_{n\geqslant 0}l^S_{\underbrace{1,1,\cdots,1}_{n}}(\varphi)u^{n}$.
Then the above equality can be read as
$$
\sum_{k\geqslant 0} (-1)^kl^S_{k+1}(\varphi)u^k=
\frac{d}{du}\log f(u).
$$
Integrating and adjusting constant terms gives
$$
\sum_{n\geqslant 0}l^S_{\underbrace{1,1,\cdots,1}_{n}}(\varphi)u^{n}
=\exp\left\{
-\sum_{n\geqslant 1} (-1)^nl^S_{n}(\varphi)\frac{u^n}{n}
\right\}
=\exp\left\{
-\sum_{n\geqslant 1} (-1)^nl^I_{n}(\varphi)\frac{u^n}{n}
\right\}
$$
because $l^S_n(\varphi)=l^I_n(\varphi)=l_n(\varphi)$ for $n>1$ and $l^S_1(\varphi)=l^I_1(\varphi)=-T$.
Since $l^I_{\bf a}(\varphi)=\frac{(-T)^m}{m!}$
for ${\bf a}=(\underbrace{1,1,\cdots,1}_{m})$,
we get $\l^S_{\bf a}(\varphi)={\mathbb L}\bigl(l^I_{\bf a}(\varphi)\bigr)$.

When ${\bf a}$ is of the form $({\bf a'},\underbrace{1,1,\cdots,1}_{l})$
with $\bf a'$ admissible,
the proof is given by the following induction on $l$.
By \eqref{series shuffle for bar},
$$
l_{\bold a'}^x(e^{TX_{51}}\varphi_{451}\varphi_{123})\cdot 
l_{\underbrace{1,1,\cdots,1}_{l}}^y(e^{TX_{51}}\varphi_{451}\varphi_{123})
={\underset{\sigma\in Sh^{\leqslant}(k,l)}{\sum}}
l_{\sigma(\bold a',\underbrace{(1,1,\cdots,1)}_{l})}^{\sigma(x,y)}
(e^{TX_{51}}\varphi_{451}\varphi_{123})
$$
with $k=dp({\bf a'})$.
By lemma \ref{lemma5} and lemma \ref{lemma6},
$$
l_{\bold a'}(\varphi)\cdot l^I_{\underbrace{1,1,\cdots,1}_{l}}(\varphi)
={\underset{\sigma\in Sh^{\leqslant}(k,l)}{\sum}}
l^I_{\sigma(\bold a',\underbrace{(1,1,\cdots,1)}_{l})}(\varphi).
$$
Then by our induction assumption, taking the image by the map $\mathbb L$ gives 
$$
l_{\bold a'}(\varphi)\cdot l^S_{\underbrace{1,1,\cdots,1}_{l}}(\varphi)
={\mathbb L}\bigl(l^I_{{\bf a'},\underbrace{1,1,\cdots,1}_{l}}(\varphi)\bigr)
+{\underset{\sigma\neq id\in Sh^{\leqslant}(k,l)}{\sum}}
l^S_{\sigma(\bold a',\underbrace{(1,1,\cdots,1)}_{l})}(\varphi).
$$
Since $l^S_{\bf a'}(\varphi)$ and $l^S_{1,\cdots,1}(\varphi)$ satisfy the series shuffle formula,
${\mathbb L}\bigl(l^I_{\bf a}(\varphi)\bigr)$ must be equal to
$\l^S_{\bf a}(\varphi)$.
\qed
\end{pf}

Embed $k\langle\langle Y_1,Y_2,\dots\rangle\rangle$ 
into $k\langle\langle X_0, X_1\rangle\rangle$
by sending $Y_m$ to $-X_0^{m-1}X_1$.
Then by the above proposition,
\begin{align*}
\l^S_{\bf a}(\varphi)
&
={\mathbb L}(l^I_{\bf a}(\varphi))
={\mathbb L}(l_{\bf a}(e^{TX_1}\varphi))
=l_{\bf a}\left({\mathbb L}(e^{TX_1}\pi_Y(\varphi))\right) 
=l_{\bf a}(\mbox{exp}
\left\{-\sum_{n=1}^{\infty}{l^I_n(\varphi)}
\frac{X_1^n}{n}\right\}\cdot\pi_Y(\varphi)) \\
&=l_{\bf a}(\mbox{exp}
\left\{-TY_1+\sum_{n=1}^{\infty}\frac{(-1)^n}{n}
c_{X_0^{n-1}X_1}(\varphi)Y^n_1\right\}
\cdot\pi_Y(\varphi))
=l_{\bf a}(e^{-TY_1}\varphi_*)
\end{align*}
for all $\bf a$ because $l_1(\varphi)=0$.
As for the third equality we use
$({\mathbb L}\otimes_k id)\circ (id\otimes_k l_{\bf a})=
(id\otimes_k l_{\bf a})\circ ({\mathbb L}\otimes_k id)$
on $k[T]\otimes_k k\langle\langle X_0, X_1\rangle\rangle$.
All $\l^S_{\bf a}(\varphi)$'s satisfy 
the series shuffle formulae \eqref{series shuffle for l},
so the $l_{\bf a}(e^{-TY_1}\varphi_*)$'s do also.
By putting $T=0$, we get that $l_{\bf a}(\varphi_*)$'s also satisfy the series shuffle formulae for all $\bf a$.
Therefore $\Delta_*(\varphi_*)=\varphi_*\widehat\otimes \varphi_*$.
This completes the proof of theorem \ref{almost main},
which implies theorem \ref{main}.
\qed

\section{Auxiliary lemmas}\label{auxiliary}
We prove all lemmas which are required to prove theorem \ref{main}
in the previous section.  

\begin{lem}\label{lemma3}
Let $\varphi$ be a group-like element in $k\langle\langle X_0,X_1\rangle\rangle$
with $c_{X_0}(\varphi)=c_{X_1}(\varphi)=0$.
Then
$l^x_{\bf a}(\varphi_{451}\varphi_{123})=l_{\bf a}(\varphi)$,
$l^y_{\bf a}(\varphi_{451}\varphi_{123})=l_{\bf a}(\varphi)$,
$l^{xy}_{\bf a}(\varphi_{451}\varphi_{123})=l_{\bf a}(\varphi)$ and
$l^{x,y}_{\bf a,\bf b}(\varphi_{451}\varphi_{123})=l_{\bf a\bf b}(\varphi)$
for any indices $\bf a$ and $\bf b$.
\end{lem}

\begin{pf}
Consider the map $\mathcal M_{0,5}\to\mathcal M_{0,4}$ sending 
$[(x_1,\cdots,x_5)]\mapsto [(x_1,x_2,x_3,x_5)]$.
This induces the projection $p_{4}: U\frak P_5 \twoheadrightarrow U\frak F_2$
sending  $X_{12}\mapsto X_0$, $X_{23}\mapsto X_1$
and $X_{i4}\mapsto 0$ ($i\in\bold Z/5$).
Express $l_{\bf a}$ as \eqref{Li-expression}.
Since $(p_4\otimes id)(\Exp\Omega_5)=\Exp\Omega_4(x)\in
U\frak F_2\widehat\otimes V({\mathcal M}_{0,5})$,
it induces the map
$p^*_4:V({\mathcal M}_{0,4})\to V({\mathcal M}_{0,5})$ which gives
$p^*_4([\frac{dz}{z}])=[\frac{dx}{x}]$ and
$p^*_4([\frac{dz}{1-z}])=[\frac{dx}{1-x}]$.
Hence $p_4^*(l_{\bf a})=l^x_{\bf a}$.
Then $l^x_{\bf a}(\varphi_{451}\varphi_{123})=
l_{\bf a}(p_4(\varphi_{451}\varphi_{123}))
=l_{\bf a}(\varphi)$
because $p_4(\varphi_{451})=0$ by our assumption $c_{X_1}(\varphi)=0$.

Next consider the map $\mathcal M_{0,5}\to\mathcal M_{0,4}$ sending 
$[(x_1,\cdots,x_5)]\mapsto [(x_1,x_3,x_4,x_5)]$.
This induces the projection $p_{2}: U\frak P_5 \twoheadrightarrow U\frak F_2$
sending  $X_{45}\mapsto X_0$, $X_{51}\mapsto X_1$
and $X_{i2}\mapsto 0$ ($i\in\bold Z/5$).
Since $(p_2\otimes id)(\Exp\Omega_5)=\Exp\Omega_4(y)\in
U\frak F_2\widehat\otimes V({\mathcal M}_{0,5})$,
it induces the map
$p^*_2:V({\mathcal M}_{0,4})\to V({\mathcal M}_{0,5})$ which gives
$p^*_2([\frac{dz}{z}])=[\frac{dy}{y}]$ and
$p^*_2([\frac{dz}{1-z}])=[\frac{dy}{1-y}]$.
Hence $p_2^*(l_{\bf a})=l^y_{\bf a}$.
Then $l^y_{\bf a}(\varphi_{451}\varphi_{123})=
l_{\bf a}(p_2(\varphi_{451}\varphi_{123}))
=l_{\bf a}(\varphi)$
because $p_2(\varphi_{123})=0$.

Similarly consider the map $\mathcal M_{0,5}\to\mathcal M_{0,4}$ sending 
$[(x_1,\cdots,x_5)]\mapsto [(x_1,x_2,x_4,x_5)]$.
This induces the projection $p_{3}: U\frak P_5 \twoheadrightarrow U\frak F_2$
sending  $X_{12}\mapsto X_0$, $X_{24}\mapsto X_1$
and $X_{i3}\mapsto 0$ ($i\in\bold Z/5$).
Since $(p_3\otimes id)(\Exp\Omega_5)=\Exp\Omega_4(xy)\in
U\frak F_2\widehat\otimes V({\mathcal M}_{0,5})$,
it induces the map
$p^*_3:V({\mathcal M}_{0,4})\to V({\mathcal M}_{0,5})$ which gives
$p^*_3([\frac{dz}{z}])=[\frac{dx}{x}+\frac{dy}{y}]$ and
$p^*_3([\frac{dz}{1-z}])=[\frac{xdy+ydx}{1-xy}]$.
Hence $p_3^*(l_{\bf a})=l^{xy}_{\bf a}$.
Then $l^{xy}_{\bf a}(\varphi_{451}\varphi_{123})=
l_{\bf a}(p_3(\varphi_{451}\varphi_{123}))
=l_{\bf a}(\varphi)$
because $p_3(\varphi_{123})=0$ by our assumption $c_{X_0}(\varphi)=0$.

Consider the embedding of Hopf algebra
$i_{123}:U\frak F_2\hookrightarrow U\frak P_5$
sending $X_0\mapsto X_{12}$ and $X_1\mapsto X_{23}$.
(Geometrically it is explained by the residue map in \cite{BF}
along the divisor $\{y=0\}$.)
Since
$(i_{123}\otimes id)(\Exp\Omega_4)=\Exp\Omega_4(z;X_{12},X_{23})
\in U\frak P_5\widehat\otimes V({\mathcal M}_{0,4})$,
it induces the map
$i^*_{123}:V({\mathcal M}_{0,5})\to V({\mathcal M}_{0,4})$ 
which gives $i^*_{123}([\frac{dy}{y}])=i^*_{123}([\frac{dy}{1-y}])
=i^*_{123}([\frac{xdy+ydx}{1-xy}])=0$.
Express $l^{x,y}_{\bf a,\bf b}$ and $l^{xy}_{\bf a}$ as \eqref{expression}.
In the expression each term contains at least one
$\frac{dy}{y}$, $\frac{dy}{1-y}$ or $\frac{xdy+ydx}{1-xy}$.
Therefore we have $i_{123}^*(l^{x,y}_{\bf a,\bf b})=0$ and
$i_{123}^*(l^{xy}_{\bf a})=0$.
Thus $l^{x,y}_{\bf a,\bf b}(\varphi_{123})=0$
and $l^{xy}_{\bf a}(\varphi_{123})=0$.
Next consider the embedding of Hopf algebra 
$i_{451}:U\frak F_2\hookrightarrow U\frak P_5$
sending $X_0\mapsto X_{45}$ and $X_1\mapsto X_{51}=X_{23}+X_{24}+X_{34}$
(geometrically caused by the divisor $\{x=1\}$.)
Since 
$(i_{451}\otimes id)(\Exp\Omega_4)=\Exp\Omega_4(z;X_{45},X_{23}+X_{24}+X_{34})
\in U\frak P_5\widehat\otimes V({\mathcal M}_{0,4})$,
it induces the map
$i^*_{451}:V({\mathcal M}_{0,5})\to V({\mathcal M}_{0,4})$
which gives
$i^*_{451}([\frac{dx}{x}])=0$,
$i^*_{451}([\frac{dx}{1-x}])=[\frac{dz}{1-z}]$,
$i^*_{451}([\frac{dy}{y}])=[\frac{dz}{z}]$,
$i^*_{451}([\frac{dy}{1-y}])=[\frac{dz}{1-z}]$ and
$i^*_{451}([\frac{xdy+ydx}{1-xy}])=[\frac{dz}{1-z}]$.
By induction on weight,
$i_{451}^*(l^{x,y}_{\bf a,\bf b})=l_{\bf a\bf b}$ and
$i_{451}^*(l^{xy}_{\bf a})=l_{\bf a}$
can be deduced from the differential equations \eqref{differential equation}:
for instance if $a_k=1$, $k\neq 1$, $b_l\neq 1$ and $l\neq 1$,
by \eqref{differential equation}
$i_{451}^*(l^{x,y}_{\bf a,\bf b})=
[\frac{dz}{1-z}|i_{451}^*(l^{x,y}_{(a_1,\cdots, a_{k-1}),{\bf b}})]
-[\frac{dz}{1-z}|i_{451}^*(l^{x,y}_{(a_1,\cdots, a_{k-1}, b_1),(b_2,\cdots, b_l)})]
+[\frac{dz}{z}|i_{451}^*(l^{x,y}_{{\bf a},(b_1,\cdots, b_l-1)})]$
but by our induction assumption it is equal to 
$[\frac{dz}{z}|l_{{\bf a},(b_1,\cdots, b_l-1)}]=l_{\bf a\bf b}$.
Thus $l^{x,y}_{\bf a,\bf b}(\varphi_{451})=l_{\bf a\bf b}(\varphi)$.
Let $\delta$ be the coproduct of $V(\mathcal M_{0,5})$.
Express $\delta(l^{x,y}_{\bf a,\bf b})=\sum_i l'_i\otimes l''_i$
with $l'_i\in V_{m'_i}$ and $l''_i\in V_{m''_i}$ for some $m'_i$ and $m''_i$.
If $m''_i\neq 0$, $l''_i(\varphi_{123})=0$
because $l''_i$ is a combination of elements
of the form $l^{x,y}_{\bf c,\bf d}$ and $l^{xy}_{\bf e}$
for some indices $\bf c$, $\bf d$ and $\bf e$.
Since $\delta(l^{x,y}_{\bf a,\bf b})(1\otimes\varphi_{451}\varphi_{123})=
\delta(l^{x,y}_{\bf a,\bf b})(\varphi_{451}\otimes\varphi_{123})$,
$l^{x,y}_{\bf a,\bf b}(\varphi_{451}\varphi_{123})=
\sum_i l'_i(\varphi_{451})\otimes l''_i(\varphi_{123})
=l^{x,y}_{\bf a,\bf b}(\varphi_{451})=l_{\bf a\bf b}(\varphi)$.
\qed
\end{pf}

\begin{lem}\label{lemma4}
Let $\varphi$ be an element as in theorem \ref{almost main}.
Suppose that $\bf a$ and $\bf b$ are admissible. Then
$l^{y,x}_{\bf a,\bf b}(\varphi_{451}\varphi_{123})=l_{\bf a\bf b}(\varphi)$.
\end{lem}

\begin{pf}
Because \eqref{5-cycle} implies $\varphi(X_0,X_1)\varphi(X_1,X_0)=1$
(take a projection $U\frak P_5 \twoheadrightarrow U\frak F_2$),
we have $\varphi_{451}\varphi_{123}=\varphi_{432}\varphi_{215}\varphi_{543}$.
Put $\delta(l^{y,x}_{\bf a,\bf b})=\sum_i l'_i\otimes l''_i$.
By the same arguments to the last paragraph of the proof of lemma \ref{lemma3},
$l^{y,x}_{\bf c,\bf d}(\varphi_{543})=0$ and $l^{xy}_{\bf e}(\varphi_{543})=0$
for any indices $\bf c$, $\bf d$ and $\bf e$.
So we have $l^{y,x}_{\bf a,\bf b}(\varphi_{451}\varphi_{123})
=l^{y,x}_{\bf a,\bf b}(\varphi_{432}\varphi_{215}\varphi_{543})
=l^{y,x}_{\bf a,\bf b}(\varphi_{432}\varphi_{215})
=\sum_i l'_i(\varphi_{432})\otimes l''_i(\varphi_{215})$.
Consider the embedding of Hopf algebra
$i_{432}:U\frak F_2\hookrightarrow U\frak P_5$
sending $X_0\mapsto X_{43}$ and $X_1\mapsto X_{32}$
(geometrically caused by the exceptional divisor obtained by blowing up at $(x,y)=(1,1)$.)
Since $(i_{432}\otimes id)(\Exp\Omega_4)=\Exp\Omega_4(z;X_{34},X_{23})
\in U\frak P_5\widehat\otimes V({\mathcal M}_{0,4})$,
it induces the morphism $i^*_{432}:V({\mathcal M}_{0,5})\to V({\mathcal M}_{0,4})$ 
which gives 
$i^*_{432}([\frac{dx}{x}])=0$,
$i^*_{432}([\frac{dx}{x-1}])=[\frac{dz}{z-1}]$,
$i^*_{432}([\frac{dy}{y}])=0$,
$i^*_{432}([\frac{dy}{y-1}])=[\frac{dz}{z}]$ and
$i^*_{432}([\frac{xdy+ydx}{1-xy}])=0$.
In each term of the expression 
$l^{y,x}_{\bf a,\bf b}
=\sum_{I=(i_m,\cdots,i_1)}c_I[\omega_{i_m}|\cdots|\omega_{i_1}]$,
the first component $\omega_{i_m}$ is always $\frac{dx}{x}$ or $\frac{dy}{y}$
because $\bf a$ and $\bf b$ are admissible.
So $i_{432}^*(l'_i)=0$ unless $m'_i=0$.
Therefore 
$\sum_i l'_i(\varphi_{432})\otimes l''_i(\varphi_{215})
=l^{y,x}_{\bf a,\bf b}(\varphi_{215})$.
By the same arguments to lemma \ref{lemma3} again,
$i_{215}^*(l^{y,x}_{\bf a,\bf b})=l_{\bf a\bf b}$.
Thus $l^{y,x}_{\bf a,\bf b}(\varphi_{215})=l_{\bf a\bf b}(\varphi)$.
\qed
\end{pf}

\begin{lem}\label{lemma5}
Let $\varphi$ be a group-like element in $k\langle\langle X_0,X_1\rangle\rangle$
with $c_{X_0}(\varphi)=c_{X_1}(\varphi)=0$.
Then 
$l^x_{\bf a}(e^{TX_{51}}\varphi_{451}\varphi_{123})=l^I_{\bf a}(\varphi)$,
$l^y_{\bf a}(e^{TX_{51}}\varphi_{451}\varphi_{123})=l^I_{\bf a}(\varphi)$,
$l^{xy}_{\bf a}(e^{TX_{51}}\varphi_{451}\varphi_{123})=l^I_{\bf a}(\varphi)$ and
$l^{x,y}_{\bf a,\bf b}(e^{TX_{51}}\varphi_{451}\varphi_{123})=l^I_{\bf a\bf b}(\varphi)$
for any index $\bf a$ and $\bf b$.
\end{lem}

\begin{pf}
By the arguments in lemma \ref{lemma3},
$$
l^x_{\bf a}(e^{TX_{51}}\varphi_{451}\varphi_{123})=
l_{\bf a}(p_4(e^{TX_{51}}\varphi_{451}\varphi_{123}))
=l_{\bf a}(e^{TX_1}\varphi)
=l^I_{\bf a}(\varphi),
$$
$$
l^y_{\bf a}(e^{TX_{51}}\varphi_{451}\varphi_{123})
=l_{\bf a}(p_2(e^{TX_{51}}\varphi_{451}\varphi_{123}))
=l_{\bf a}(e^{TX_1}\varphi)
=l^I_{\bf a}(\varphi),
$$
$$
l^{xy}_{\bf a}(e^{TX_{51}}\varphi_{451}\varphi_{123})
=l_{\bf a}(p_3(e^{TX_{51}}\varphi_{451}\varphi_{123}))
=l_{\bf a}(e^{TX_1}\varphi)
=l^I_{\bf a}(\varphi)
$$
\qquad and
$
l^{x,y}_{\bf a,\bf b}(e^{TX_{51}}\varphi_{451}\varphi_{123})
=l^{x,y}_{\bf a,\bf b}(e^{TX_{51}}\varphi_{451})
=l_{\bf a\bf b}(e^{TX_1}\varphi)
=l^I_{\bf a\bf b}(\varphi).
$
\qed
\end{pf}

\begin{lem}\label{lemma6}
Let $\varphi$ be an element as in theorem \ref{almost main}.
Suppose that $\bf b$ is admissible. Then
$l^{y,x}_{\bf a,\bf b}(e^{TX_{51}}\varphi_{451}\varphi_{123})=l_{\bf a\bf b}(\varphi)$.
\end{lem}

\begin{pf}
We have $l^{y,x}_{\bf a,\bf b}(e^{TX_{51}}\varphi_{451}\varphi_{123})
=l^{y,x}_{\bf a,\bf b}(e^{TX_{51}}\varphi_{432}\varphi_{215}\varphi_{543})
=l^{y,x}_{\bf a,\bf b}(\varphi_{432}e^{TX_{51}}\varphi_{215})
=\sum_i l'_i(\varphi_{432})\otimes l''_i(e^{TX_{51}}\varphi_{215})$
because $e^{TX_{51}}\varphi_{432}=\varphi_{432}e^{TX_{51}}$.
Since $\bf b$ is admissible,
$i^*_{432}(l'_i)$ is of the form $\alpha[\frac{dz}{z}|\cdots|\frac{dz}{z}]$
with $\alpha\in k$.
By our assumption, $c_{X_0\cdots X_0}(\varphi)=0$. 
So $l'_i(i_{432}(\varphi))=l'_i(\varphi_{432})=0$
unless $m'_i=0$.
Thus 
$\sum_i l'_i(\varphi_{432})\otimes l''_i(e^{TX_{51}}\varphi_{215})
=l^{y,x}_{\bf a,\bf b}(e^{TX_{51}}\varphi_{215})
=l_{\bf a\bf b}(e^{TX_1}\varphi)$.
Since $\bf b$ is admissible,
$l_{\bf a\bf b}(e^{TX_1}\varphi)=l_{\bf a\bf b}(\varphi)$.
\qed
\end{pf}

\section{Proof of theorem \ref{gamma factorization}}
\label{proof of gamma factorization}
This section gives a proof of theorem \ref{gamma factorization}.

Put $F'_2(k)=\{\varphi\in U\frak F_2| \varphi:\text{group-like}, 
c_{X_0}(\varphi)=c_{X_1}(\varphi)=0\}$.
It forms a group with respect to \eqref{multiplication}
and contains $DMR_0(k)$ as a subgroup.
Consider the map $m:(F'_2(k),\circ)\to k[[x_0,x_1]]^\times$
sending $\varphi$ to its meta-abelian quotient $B_\varphi(x_0,x_1)$.
By a direct calculation we see that it is a group homomorphism, i.e.
$m(\varphi_1\circ\varphi_2)=m(\varphi_1)\cdot m(\varphi_2)$.
Put $B(k)=\bigl\{b\in k[[x_0,x_1]]^\times\bigl|
b(x_0,x_1)=\frac{c(x_0)c(x_1)}{c(x_0+x_1)}
\text{ for } c(s)\in 1+s^2 k[[s]]\bigr\}$.
It is a subgroup of $k[[x_0,x_1]]^\times$.
The first statement of theorem \ref{gamma factorization}
claims $m(DMR_\mu(k))\subset B(k)$ for all $\mu\in k$.

\begin{prop}\label{0-case}
$m(DMR_0(k))\subset B(k)$.
\end{prop}

\begin{pf}
Let $M$ be the Lie algebra homomorphism,
associated with $m|_{DMR_0}$,
from the Lie algebra $\frak{dmr}_0$ (see Appendix \ref{appendix}) of $DMR_0$ to the trivial Lie 
algebra $k[[x_0,x_1]]$.
In order to prove our proposition
it is enough to show $M(\frak{dmr}_0)\subset\frak B$
where $\frak B$ is the Lie subalgebra
$\bigl\{\beta\in k[[x_0,x_1]]\bigl|
\beta(x_0,x_1)=\gamma(x_0)+\gamma(x_1)-\gamma(x_0+x_1)
\text{ for } \gamma(s)\in s^2 k[[s]]\bigr\}$
with trivial Lie structure. 
\end{pf}

\begin{lem}
For $\psi\in \frak{dmr}_0$ with $\psi=\psi_{X_0}X_0+\psi_{X_1}X_1$, $M(\psi)=(\psi_{X_1}X_1)^\text{ab}$.
\end{lem}

\begin{pf}
The exponential map $\Exp:\frak{dmr}_0\to DMR_0$ is given by the formula
$\psi\mapsto\sum_{i=0}^\infty\frac{1}{i!}(\mu_\psi+d_\psi)^i(1)
=1+\psi+\frac{1}{2}(\psi^2+d_\psi(\psi))^2+\frac{1}{6}
(\psi^3+2\psi d_\psi(\psi)+d_\psi(\psi)\psi+d_\psi^2(\psi))
+\cdots$
($\mu_\psi$: the left multiplication by $\psi$
and $d_\psi$: see our appendix)
in \cite{DG} Remark 5.15.
By a direct calculation, it can be checked that
$M\circ\Exp(\psi)=\exp(\psi_{X_1}X_1)^\text{ab}$,
which implies our lemma.
\qed
\end{pf}

By the above lemma the proof of $M(\frak{dmr}_0)\subset\frak B$
is reduced to the following.

\begin{lem}
For $\psi\in \frak{dmr}_0$,
$$
\underset{dp({\bf a})=m}{\sum_{wt({\bf a})=w}}l_{\bf a}(\psi)=
\begin{cases}
(-1)^{m-1}\binom{w}{m}\frac{l_w(\psi)}{w} &\text{ for }m<w,\\
0 &\text{ for }m=w.
\end{cases}
$$
\end{lem}

\begin{pf}
By \eqref{Lie double shuffle} we have
$\sum_{\sigma\in Sh^{\leqslant}(k,l)}l_{\sigma(\bold a,\bold b)}(\psi_*)=0$
with $dp({\bf a})=k$ and $dp({\bf b})=l$.
Summing up all pairs $(\bold a,\bold b)$ with
$wt({\bf a})=k$, $dp({\bf a})=1$ and $wt({\bf a})+wt({\bf b})=w$,
we get
$$
\sum_{wt({\bf a})=w,dp({\bf a})=k+1}(k+1)l_{\bf a}(\psi_*)+
\sum_{wt({\bf a})=w,dp({\bf a})=k}(w-k)l_{\bf a}(\psi_*)=0.
$$
Then the lemma follows by induction because $l_{\bf a}(\psi_*)=l_{\bf a}(\psi)$
for $dp({\bf a})\neq wt({\bf a})$.
\qed
\end{pf}

The above lemma implies proposition \ref{0-case}.\qed

To complete the proof of the first statement of theorem \ref{gamma factorization}
we may assume that $\mu=1$ and $k=\bf Q$.
Let $\varphi$ be any element of $DMR_1(\bf Q)$.
Put $\varphi_{KZ}=\varPhi_{KZ}(\frac{X_0}{2\pi i}, \frac{X_1}{2\pi i})$.
It gives an element of $DMR_1(\bf C)$.
Because $DMR_1$ is a right $DMR_0$-torsor
there exists uniquely $\varphi'\in DMR_0(\bf C)$
such that $\varphi=\varphi_{KZ}\circ\varphi'$.
In \cite{Dr} it is shown that $m(\varphi_{KZ})\in B(\bf C)$.
By proposition \ref{0-case}, $m(\varphi')\in B(\bf C)$.
Thus $m(\varphi)=m(\varphi_{KZ})m(\varphi')$ must lie on $B$.
The proof of the second statement is easy. 
Express $\log\Gamma_\varphi(s)=\sum_{n=2}^\infty d_n(\varphi)s^n$.
Then by \eqref{gamma factorization formula}
$d_n(\varphi)=\frac{-1}{n}c_{X_0^{n-1}X_1}(\varphi)$.
This completes the proof of theorem \ref{gamma factorization}.
\qed

\appendix
\section{Review of the proof of Racinet's theorem}\label{appendix}
In \cite{R} theorem I, Racinet shows a highly non-trivial result
that $DMR_0$ is closed under the multiplication \eqref{multiplication}.
However his proof looks too complicated. 
The aim of this appendix is to review the essential part (\cite{R} proposition 4.A.i)
of his proof clearly
in the case of $\Gamma=\{1\}$
in order to help the readers to catch his arguments.

In \cite{R}3.3.1, $\frak{dmr}_0$ is introduced to be the set of formal Lie series 
$\psi\in \frak F_2=\{\psi\in U\frak F_2|\Delta(\psi)=1\otimes\psi+\psi\otimes 1\}$
satisfying $c_{X_0}(\psi)=c_{X_1}(\psi)=0$ and
\begin{equation}\label{Lie double shuffle}
\Delta_*(\psi_*)=1\otimes\psi_*+\psi_*\otimes 1
\end{equation}
with
$\psi_*=\psi_\text{corr}+\pi_Y(\psi)$ and
$\psi_{\text{corr}}=\sum_{n=1}^{\infty}
\frac{(-1)^n}{n}c_{X_0^{n-1}X_1}(\psi)Y_1^n$.
It is the tangent vector space at the origin of $DMR_0$.

\begin{thm}[\cite{R} proposition 4.A.i]\label{4.A.i}
The set $\frak{dmr}_0$ has a structure of Lie algebra with 
the Lie bracket
\footnote{
For our convenience,
we change the order of bracket in \cite{R} (3.1.10.2.)
}
given by 
$$
\{\psi_1,\psi_2\}=d_{\psi_2}(\psi_1)-d_{\psi_1}(\psi_2)-[\psi_1,\psi_2]
$$
where $d_\psi$ ($\psi\in U\frak F_2$) is the derivation of $U\frak F_2$ given by
$d_\psi(X_0)=0$ and $d_\psi(X_1)=[X_1,\psi]$.
\end{thm}

\begin{pf}
Put $U \frak F_Y=k\langle\langle Y_1,Y_2,\dots\rangle\rangle$.
It is the universal enveloping algebra of the Lie algebra
$\frak F_Y=\{\psi_*\in U\frak F_Y|
\psi_* \text{ satisfies }\eqref{Lie double shuffle}\}$
with the coproduct $\Delta_*$.
By the algebraic map sending $Y_m$ to $-X_0^{m-1}X_1$ ($m=1,2,\dots$),
we frequently regard $U\frak F_Y$ as a $k$-linear (or algebraic) subspace of $U\frak F_2$.
Define $S_X$ to be the antipode, the {\it anti}-automorphism
of $U\frak F_2$ such that
$S_X(X_0)=-X_0$ and $S_X(X_1)=-X_1$ 
and $R_Y$ to be the {\it anti}-automorphism of $U\frak F_Y$ such that $R_Y(Y_n)=Y_n$
($n\geqslant 1$).
It is easily proved that $S_X(x_1w)=(-1)^{p+1}R_Y(x_1w)$
for $w\in U\frak F_Y$ with degree $p$ and
that $S_X(f)=-f$ (resp. $R_Y(f)=-f$)
for $f\in\frak F_2$ (resp. $f\in\frak F_Y$).

For a homogeneous element of $U\frak F_2$ with degree $p$,
we denote $\{f_i\}^p_{i=0}$ to be elements of $U\frak F_Y$ characterised by
$f=\sum^p_{i=0}f_iX_0^i$ and
define $f_{i,j}=f_iY_j+Y_j\bar f_i\in U\frak F_2$
($0\leqslant i\leqslant p$, $0\leqslant j$)
with $\bar f_i=(-1)^p R_Y(f_i)$.

Define the $k$-linear endomorphism $s_f$ 
($f\in U\frak F_2$)
of $U\frak F_2$ by $s_f(v)=fv+d_f(v)$.
It induces the $k$-linear endomorphism $s^Y_f$ of $U\frak F_Y$ such that
$\pi_Y\circ s_f=s^Y_f\circ\pi_Y$.
Define the $k$-linear endomorphism $D^Y_f$ of $U\frak F_Y$ 
by $D^Y_f(w)=s^Y_f(w)-w\pi_Y(f)$.
Direct computations show that $D^Y_f$ forms a derivation of
the non-commutative algebra $U\frak F_Y$ 
and actually it is the restriction into $U\frak F_Y$ of the derivation $D_f$
of $U\frak F_2$ defined by $X_0\mapsto [f,X_0]$ and $X_1\mapsto 0$.

\begin{lem}
\label{4.3.4}
Let $f$ be a homogeneous element of $U\frak F_2$ with degree $p(\geqslant 1)$
such that
$S_X(f)=-f$.
Then for $n\geqslant 1$
$$D^Y_f(Y_n)=X_0^{n-1}fX_1-fX_0^{n-1}X_1
=\sum_{i=0}^p f_{i,i+n}.$$
\end{lem}

\begin{pf}
The first equality is easy.
Denote $f_i=\sum_{j=0}^{p-1-i}X_0^jX_1f^j_i$ ($0\leqslant i\leqslant p-1$).
By $S_X(f)=-f$, we have
$X_1f^j_i=(-1)^{p+1}R_Y(X_1f^i_j)$
($0\leqslant i\leqslant p-1$, $0\leqslant j\leqslant p-1-i$)
and $\{1+(-1)^p\}f_p=0$.
Hence
$\bar f_i=-\sum_{j=0}^{p-1-i}f^i_jX_0^jX_1$ ($0\leqslant i\leqslant p-1$) and
$\bar f_p=(-1)^pf_p$.
Then direct calculation using them yields the second equality.
\qed
\end{pf}

Put $\partial_0$ to be the derivation of $U\frak F_2$ defined by $\partial_0(X_0)=1$
and $\partial_0(X_1)=0$.
Define the map $\sec:U\frak F_Y\to U\frak F_2$ by
$\sec(w)=\sum_{i\geqslant 0}\frac{(-1)^i}{i!}\partial^i_0(w)X_0^i$,
which actually maps to  $\ker\partial_0\subset U\frak F_2$.
It is easy to see that
the map $\sec:U\frak F_2\to\ker\partial_0$ is the inverse of $\pi_Y|_{\ker\partial_0}$. 

\begin{lem}\label{4.3.9}
Let $g$ be a homogeneous element of $U\frak F_2$ with degree $p(\geqslant 1)$.
Assume $g\in \frak F_Y$. Set $f=\sec g$.
Then for $n\geqslant 1$
$$\Delta_*\circ D^Y_f(Y_n)-(id\otimes D^Y_f+D^Y_f\otimes id)\Delta_*(Y_n)=
\sum_{k=0}^p\sum_{i=k}^p\{f_{i,i-k}\otimes Y_{n+k}
+Y_{n+k}\otimes f_{i,i-k}\}.$$
\end{lem}

\begin{pf}
It is well-known that if for $i\geqslant 1$ we define $U_i$
to be the weight $i$ part of $\log(1+Y_1+Y_2+\cdots)$
(so $U_1=Y_1$, $U_2=Y_2-\frac{Y_1^2}{2}$, etc),
then $g\in U\frak F_Y$ lies on $\frak F_Y$
if and only if
$g$ lies on the free Lie algebra generated by $U_1,U_2,\dots$.
By direct calculation, $\partial_0(U_n)=(n-1)U_{n-1}$ ($n\geqslant 1$).
So by induction, we get $\partial_0(w)\in\frak F_Y$ if $w\in \frak F_Y$.
Because $f_i=\frac{(-1)^i}{i!}\partial^i_0(g)$ and
$\bar f_i=(-1)^p R_Y(f_i)$ ($0\leqslant i \leqslant p$),
both $f_i$ and $\bar f_i$ satisfy \eqref{Lie double shuffle}
(belong to $\frak F_Y$).
Then it follows that
(the first term)$=\sum_{i=0}^p\Delta_*(f_iY_{i+n}+Y_{i+n}\bar f_i)
=\sum_{i=0}^p\sum_{k=0}^{n+i}
(f_{i,i+n-k}\otimes Y_{k}+Y_{k}\otimes f_{i,i+n-k})$
and
(the second term)$=\sum_{k=0}^{n-1}\sum_{i=0}^p
(f_{i,i+n-k}\otimes Y_{k}+Y_{k}\otimes f_{i,i+n-k}).
$
Therefore
(l.h.s)$=\sum_{i=0}^p\sum_{k=0}^{i}
(f_{i,i-k}\otimes Y_{n+k}+Y_{n+k}\otimes f_{i,i-k})
=\sum_{k=0}^p\sum_{i=k}^{p}
(f_{i,i-k}\otimes Y_{n+k}+Y_{n+k}\otimes f_{i,i-k}).
$
\qed
\end{pf}

\begin{lem}\label{4.3.8}
Let $g$ be a homogeneous element of $U\frak F_2$ with degree $p(\geqslant 1)$.
Set $f=\sec g$. 
Assume $g\in \frak F_Y$ and $S_X(f)=-f$.
Then

$\sum_{i=k}^pf_{i,i-k}=
\begin{cases}
(-1)^{p-k-1}\binom{p-1}{k}\{1+(-1)^p\}c_{X_0^{p-1}X_1}(g)Y_{p-k}
\qquad (0\leqslant k \leqslant p-1),\\
0\qquad\qquad\qquad\qquad\qquad\qquad\qquad\qquad\qquad\qquad (k=p).\\
\end{cases}
$
\end{lem}

\begin{pf}
Let $\partial_{U_i}$ ($i\geqslant 1$) be the derivation defined by
$1$ on $U_i$ and $0$ on $U_j$ ($j\neq i$).
It is not difficult to show that $\partial_{U_i}(w)=c_{Y_i}(w)$ 
for $w\in\frak F_Y$ (cf. \cite{R} proposition 2.3.8).
By $g\in\frak F_Y$,  we  have $f_i,\bar f_i\in\frak F_Y$. 
Whence direct calculation yields for $k\geqslant 0$
$$
\partial_{U_{k+1}}(\sum_{i=0}^p f_{i,i+1})=
\sum_{i=0}^pc_{Y_{k+1}}(f_i+\bar f_i)Y_{i+1}+
\sum_{i=k}^pf_{i,i-k}.
$$
By $\deg f_i=p-i$, $\bar f_i=(-1)^p R_Y(f_i)$ and
$f_i=\frac{(-1)^i}{i!}\partial^i_0(g)$,
$\sum_{i=0}^pc_{Y_{k+1}}(f_i+\bar f_i)Y_{i+1}
=-\{1+(-1)^p\}c_{X_0^kX_1}(f_{p-k-1})Y_{p-k}
=\{1+(-1)^p\}\frac{(-1)^{p-k}}{(p-k-1)!}
\frac{(p-1)!}{k!}c_{X_0^{p-1}X_1}(f)Y_{p-k}
$
for $0\leqslant k\leqslant p-1$.
By the definitions above it is immediate that $D_f(Y_1)=0$.
So by lemma \ref{4.3.4}, $\sum_{i=0}^pf_{i,i+1}=0$.
It implies the desired equality for $0\leqslant k \leqslant p-1$.
The case for $k=p$ is immediate.
\qed
\end{pf}

\begin{prop}[\cite{R} proposition 4.3.1]\label{4.3.1}
For $\psi\in\frak{dmr}_0$, $s^Y_{\sec\psi_*}$ forms a coderivation
of $U\frak F_Y$ with respect to $\Delta_*$.
\end{prop}

\begin{pf}
Because
$\frak F_2^{>1}\subset\ker\partial_0$ and $\sec(Y^n_1)=Y^n_1$ ($n\geqslant 1$),
we have
\begin{equation}\label{section}
\sec\psi_*=\psi+\psi_\text{corr}
\end{equation}
for $\psi\in\frak{dmr}_0$.
Since $\psi$ is an element of $\frak F_2$, $S_X(\psi)=-\psi$.
Because $\psi$ lies in $\frak{dmr}_0$, it is known that
$c_{X_0^{n-1}X_1}(\psi)=0$ for even $n$ (cf. \cite{R} proposition 3.3.3).
So $S_X(\psi_\text{corr})=-\psi_\text{corr}$.
Therefore, $S_X(\sec\psi_*)=-\sec\psi_*$ by \eqref{section}
and we can apply lemma \ref{4.3.9} and lemma \ref{4.3.8}
with $g=\psi_*$ and $f=\sec\psi_*$
to obtain the identity
$\Delta_*\circ D^Y_{\sec\psi_*}(Y_n)=(id\otimes D^Y_{\sec\psi_*}+D^Y_{\sec\psi_*}\otimes id)\Delta_*(Y_n)$
for all $n\geqslant 1$.
This implies that the derivation $D^Y_{\sec\psi_*}$ forms a coderivation.
Since $\psi_*=\pi_Y\circ\sec\psi_*$ lies on $\frak F_Y$,
the right multiplication by $\psi_*$ forms a coderivation.
Therefore  $s^Y_{\sec\psi_*}$ must form a coderivation.
\qed
\end{pf}

Recall that $s^Y_\psi(1)=\pi_Y(\psi)$ and 
$[s_\psi,s_{X_1^n}]=0$ for any $\psi$ and $n>0$
and
$s_{\{\psi_1,\psi_2\}}=[s_{\psi_2},s_{\psi_1}]$ for any $\psi_1$ and $\psi_2$.
Assume that $\psi_1$ and $\psi_2$ be elements of $\frak{dmr}_0$.
Then
$s_{\{\psi_1,\psi_2\}}=[s_{\psi_2},s_{\psi_1}] 
=[s_{\sec\psi_{2*}},s_{\sec\psi_{1*}}]$ by \eqref{section}.
Therefore
$\pi_Y(\{\psi_1,\psi_2\})=s^Y_{\{\psi_1,\psi_2\}}(1)
=[s^Y_{\sec\psi_{2*}},s^Y_{\sec\psi_{1*}}](1)$.
By proposition \ref{4.3.1}, $s^Y_{\sec\psi_{i*}}$ ($i=1,2$) forms a coderivation.
So $[s^Y_{\sec\psi_{2*}},s^Y_{\sec\psi_{1*}}]$ forms a coderivation.
Since $\pi_Y(\{\psi_1,\psi_2\})$  is the image of $1$ 
by the coderivation, it must be Lie-like with respect to $\Delta_*$.
That means
\eqref{Lie double shuffle} holds for $\psi=\{\psi_1,\psi_2\}$
because $\psi_\text{corr}=0$. 
It completes the proof of theorem \ref{4.A.i}.
\qed \end{pf}

\begin{ack}The author is grateful to F.Brown for valuable discussions and P.Deligne for his suggestions. He also has  special thanks to L.Schneps for her significant comments and her explaining  me a shorter proof in the Lie algebra setting. He is supported by JSPS Postdoctoral Fellowships for Research Abroad and JSPS Core-to-Core Program 18005.\end{ack}


\end{document}